\documentclass[12pt]{article}
\usepackage[cp1251]{inputenc}
\usepackage[russian]{babel}
\usepackage{amssymb}
\usepackage{amsmath}

\begin{document}

UDC 517.51


\begin{center}
\textbf{\large A short and simple proof of the Jurkat--Waterman
theorem on conjugate functions}
\end{center}
\begin{center}
V. V. Lebedev
\end{center}

\begin{quotation}
{\small \textsc{Abstract.} It is well--known that certain
properties of continuous functions on the circle $\mathbb T$,
related to the Fourier expansion, can be improved by a change
of variable, i.e., by a homeomorphism of the circle onto
itself. One of the results in this area is the Jurkat--Waterman
theorem on conjugate functions, which improves the classical
Bohr--P\'al theorem.  In the present work we provide a short
and technically very simple proof of the Jurkat--Waterman
theorem. Our approach yields a stronger result.

References: 14 items.

}
\end{quotation}

\quad

It is well--known that certain properties of continuous
functions $f$ on the circle $\mathbb T=\mathbb R/2\pi\mathbb
Z$, related to their Fourier expansions
$$
f(t)\sim\sum_{k\in\mathbb Z}\widehat{f}(k)e^{ikt},
$$
can be improved by a change of variable, i.e., by a
homeomorphism of the circle onto itself. The first result in
this area is the classical Bohr--P\'al theorem [1], [2], which
states that  for every real-valued continuous function $f$ on
$\mathbb T$ there exists a homeomorphism $h$ of $\mathbb T$
onto itself such that the superposition $f\circ h$ belongs to
the space $U(\mathbb T)$ of functions with uniformly convergent
Fourier series. \footnote{This theorem was obtained by J. P\'al
in a somewhat weaker form, the final version is due to H. Bohr.
The proof was subsequently simplified by R. Salem in [3], his
proof can also be found in [4, Ch. IV, Sec. 12].}

W. Jurkat and D. Waterman [5] (see also [6, Theorem 9.5])
obtained the following improvement of the Bohr--P\'al theorem.

\quad

\textsc{Theorem.} \emph{Let $f$ be a real-valued continuous
function on the circle $\mathbb T$. Then there exists a
homeomorphism $h$ of $\mathbb T$ onto itself such that the
function $\widetilde{f\circ h}$ conjugate to the superposition
$f\circ h$ is continuous and is of bounded variation on
$\mathbb T$.}

\quad

We recall that for a function $g\in L^1(\mathbb T)$ the
conjugate function $\widetilde{g}$ is defined by
$$
\widetilde{g}(t)=\lim_{\varepsilon\rightarrow +0}
\frac{1}{\pi}\int_{\varepsilon\leq |t-\theta|\leq\pi}\frac{g(\theta)}{2\tg\frac{t-\theta}{2}}d\theta,
\qquad t\in\mathbb T
$$
(the finite limit exists for almost all $t\in\mathbb T$). If
$g\in L^2(\mathbb T)$, then $\widetilde{g}\in L^2(\mathbb T)$,
and the Fourier coefficients of the functions $g,
\widetilde{g}$ are related by
$\widehat{\widetilde{g}}(k)=-i(\mathrm{sgn}\,k)\widehat{g}(k)$
(here $\mathrm{sgn}\,k=k/|k|$ for $k\neq 0$ and
$\mathrm{sgn}\,0=0$). Consider functions analytic in the disk
$D=\{z\in \mathbb C : |z|<1\}$ of the complex plane $\mathbb
C$. Let $F$ be a function of the Hardy class $H^2(D)$, then,
setting $u(t)=\mathrm{Re}\,F(e^{it}),
v(t)=\mathrm{Im}\,F(e^{it})$, we have
$\widetilde{u}=v-\mathrm{Im}\,F(0)$.

It is well-known that every continuous function of bounded
variation is in $U(\mathbb T)$. It is also known (see, e.g.,
[4, Ch. VIII, Sec. 19]) that if both $g$ and $\widetilde{g}$
are continuous and one of them is in $U(\mathbb T)$, then the
other one is in $U(\mathbb T)$ as well; thus, the
Jurkat--Waterman theorem implies the Bohr--P\'al theorem,
moreover we obtain $f\circ h\in U(\mathbb T), \widetilde{f\circ
h}\in U(\mathbb T)$.

In the present paper we give an alternative proof of the
Jurkat--Waterman theorem. Our proof is based on general facts
and, being void of technical details, is extremely short and
simple. Similarly to the original proof of the Bohr--P\'al
theorem, and similarly to the original proof of the
Jurkat--Waterman theorem, our proof also involves a conformal
mapping of the disk $D$ onto a certain domain in $\mathbb C$. A
simplification of the proof is achieved due to a special choice
of the domain, which is different from the one used in the
original works. In addition, combining our method with a result
of the work [7], we obtain a stronger assertion, see Remark 1
after the end of the proof of the theorem.

\quad

\textsc{Proof of the Theorem.} Consider a curve
$\gamma\subseteq\mathbb C$ given by $\gamma(t)=e^{f(t)+it},
t\in [0, 2\pi]$. This is a closed curve without
self-intersections. Let $V$ be the inner domain of this curve.
We have $0\in V$. Consider a conformal mapping $G$ of the disk
$D$ onto $V$. We can assume that $G(0)=0$. As is known the
mapping $G$ extends to a homeomorphism of the closure
$\overline{D}$ of $D$ onto the closure $\overline{V}$ of $V$
and, being thus extended, homeomorphically maps the boundary
$\partial D$ of $D$ onto the boundary $\gamma$ of $V$.
Preserving our notation we denote this extension by $G$. It is
clear that there exists a homeomorphism $h$ of the segment $[0,
2\pi]$ onto itself such that $G(e^{it})=\gamma(h(t))$. Thus,
$$
G(e^{it})=e^{f\circ h(t)+ih(t)}, \qquad t\in [0, 2\pi].
\eqno(1)
$$

Since $G(0)=0$, we have $G(z)=z\xi(z), z\in D$, where $\xi$ is
a function analytic in $D$. It is clear that $\xi$ extends to a
function continuous in $\overline{D}$. Preserving our notation
we denote this extension by $\xi$. Note now that if $z\in
\overline{D}\setminus\{0\}$, then $\xi(z)\neq 0$, for otherwise
we would have a point $z_0\in\overline{D}, z_0\neq 0,$ such
that $G(z_0)=0$, which is impossible since $G(0)=0$ and $G$ is
one-to-one. At the origin we have $\xi(0)=G'(0)\neq 0$, since
$G$ is conformal. Thus, the function $\xi$ is analytic in $D$,
continuous in $\overline{D}$, and non-vanishing in
$\overline{D}$. It follows that $\xi$ is of the form
$\xi(z)=e^{\eta(z)}, z\in\overline{D}$, where $\eta$ is a
function analytic in $D$ and continuous in $\overline{D}$.
Indeed, one can take $\eta$ to be any fixed continuous branch
of $\log\xi(z), z\in D$. Since $\xi$ is continuous in
$\overline{D}$ and does not vanish in $\overline{D}$, the
function $\eta$ so defined is continuous in $\overline{D}$ as
well. \footnote{For another proof of the existence of the
representation $\xi=e^\eta$, based on the theory of commutative
Banach algebras, see [8, Theorem 2.11].}

Taking (1) into account, we have
$$
e^{\eta(e^{it})}=\xi(e^{it})=e^{-it}G(e^{it})=e^{-it+f\circ h(t)+ih(t)}.
$$
Hence,
$$
\eta(e^{it})=-it+f\circ h(t)+ih(t)+ i2\pi N(t),
\eqno(2)
$$
where $N(t)$ takes integer values for all $t$. Since the
function $\eta(e^{it})$ is continuous, the function $N(t)$ is
continuous as well and hence it is identically constant,
$N(t)=N_0$ for all $t\in [0, 2\pi]$. Thus, from (2) it follows
that
$$
\eta(e^{it})=f\circ h(t)+i(h(t)-t+ 2\pi N_0), \qquad t\in[0, 2\pi].
$$
Hence,
$$
\widetilde{f\circ h}(t)=h(t)-t+ 2\pi N_0-\mathrm{Im}\,\eta(0).
\qquad t\in[0, 2\pi].
\eqno(3)
$$
Identifying $h$ on the left-hand side of (3) with a
homeomorphism of $\mathbb T$ we complete the proof. $\Box$

\quad

\textsc{Remarks.} 1. The above method of proof of the
Jurkat--Waterman theorem yields a stronger result, namely:
\emph{for every real-valued continuous function $f$ on $\mathbb
T$ there exists a homeomorphism $h$ of $\mathbb T$ onto itself
such that $\widetilde{f\circ h}$ is of bounded variation and
has the logarithmic modulus of continuity. In addition, the
homeomorphism $h$ has the logarithmic modulus of continuity.}
As usual, saying that a function $g$ has the logarithmic
modulus of continuity we mean that
$$
\sup_{|t_1-t_2|\leq\delta}|g(t_1)-g(t_2)|=O\bigg(\frac{1}{\log1/\delta}\bigg),
\qquad \delta\rightarrow+0.
$$
Let us verify that the stated improvement holds. The domain
$V$, that we considered in the proof of the theorem, is
star-like with respect to $0$, i.e., for every point $z\in V$
the line segment $[0, z]$ lies in $V$. Taking relations (1),
(3) into account it remains to use the following result by N.
S. Stylianopoulos and E. Wegert [7]: let $V\subseteq\mathbb C$
be the inner domain of a Jordan curve, and let $G$ be a
conformal mapping of the disk $D$ onto $V$ with $G(0)=0$, then
assuming $V$ is star-like with respect to $0$, the function
$$
\varphi(t)=\arg G(e^{it}), \qquad  t\in [0, 2\pi],
$$
has the logarithmic modulus of continuity.

2. As is noted by Jurkat and Waterman [5] their theorem implies
that for every real-valued continuous function $f$ on $\mathbb
T$ there exists a homeomorphism $h$ such that the Fourier
coefficients of the superposition $f\circ h$ are rapidly
decreasing, namely $\widehat{f\circ h}(k)=O(1/|k|)$. This is a
direct consequence of a well-known estimate for the Fourier
coefficients of a function of bounded variation and of the
relation |$\widehat{g}(k)|=|\widehat{\widetilde{g}}(k)|,\,k\neq
0$.

3. Recall that the Sobolev space $W_2^{1/2}(\mathbb T)$
consists of (integrable) functions $g$ on $\mathbb T$
satisfying $\sum_{k\in\mathbb Z}|\widehat{g}(k)|^2|k|<\infty$.
It is easy to see that from the Jurkat--Waterman theorem it
follows that for every real-valued continuous function $f$ on
$\mathbb T$ there exists a homeomorphism $h$ such that $f\circ
h\in W_2^{1/2}(\mathbb T)$. Indeed, if $g$ is a real-valued
continuous function such that $\widetilde{g}$ is continuous and
of bounded variation, then $g\in W_2^{1/2}(\mathbb T)$. One can
verify this as follows. Taking into account that
$\widehat{\widetilde{g}}(k)=-i(\mathrm{sgn}\,k)\widehat{g}(k)$,
we have
$$
\frac{1}{2\pi}\int_\mathbb T e^{ikt} d\widetilde{g}(t)=
-\frac{1}{2\pi}\int_\mathbb T \widetilde{g}(t)de^{ikt}=
-ik\overline{\widehat{\widetilde{g}}(k)}=|k|\overline{\widehat{g}(k)}
\eqno(4)
$$
(the dash stands for the complex conjugation). Let
$\sigma_N(g)$ be the Fej\'er sums of $g$, i.e.,
$$
\sigma_N(g)(t)=\sum_{|k|\leq N}\widehat{g}(k)\bigg(1-\frac{|k|}{N}\bigg)e^{ikt}.
$$
Using (4), we see that
$$
\frac{1}{2\pi}\int_\mathbb T\sigma_N(g)(t)d\widetilde{g}(t)=
\sum_{|k|\leq N}|\widehat{g}(k)|^2|k|\bigg(1-\frac{|k|}{N}\bigg).
$$
Passing to the limit as $N\rightarrow\infty$, and taking into
account that $\sigma_N(g)(t)$ converge uniformly to $g(t)$, we
obtain $g\in W_2^{1/2}(\mathbb T)$.\footnote{In addition we
obtain $\frac{1}{2\pi}\int_\mathbb T g(t) d\widetilde{g}(t)=
\sum_{k\in\mathbb Z}|\widehat{g}(k)|^2|k|$.} Note by the way
that for every continuous function that belongs to
$W_2^{1/2}(\mathbb T)$, we have $g\in U(\mathbb T)$.

4. The Bohr--P\'al theorem was extended to families of
functions by J. -P. Kahane and Y. Katznelson. These authors
showed [9] (see also [10]) that if $K$ is a compact family in
the space $C(\mathbb T)$ of continuous functions on $\mathbb
T$, then there exists a homeomorphism $h$ such that $f\circ
h\in U(\mathbb T)$ for all $f\in K$. In particular the
assertion of the Bohr--P\'al theorem holds for a complex-valued
$f$. On the other hand it is impossible to get a similar
extension of the Jurkat--Waterman theorem. Moreover there exist
two real-valued continuous functions $u$ and $v$ such that
there is no homeomorphism $h$ of $\mathbb T$ onto itself with
the property that both inclusions $u\circ h\in
W_2^{1/2}(\mathbb T), v\circ h\in W_2^{1/2}(\mathbb T)$ hold
simultaneously. This was first established in [11, Theorem 4].
Thus, there exists a complex-valued continuous function $f$
such that $f\circ h\notin W_2^{1/2}(\mathbb T)$, for every
homeomorphism $h$ (see also [12] where it is shown that for
each $\alpha<1/2$ such an $f$ can be chosen so that it
satisfies the  Lipschitz condition of order $\alpha$).

5. Another improvement of the Bohr--P\'al theorem is due to A.
A. Sahakian [13]. He showed that if $\{\varepsilon(n), n=1, 2,
\ldots\}$ is a positive sequence monotonically decreasing to
zero, with $\sum_{n=1}^\infty\varepsilon(n)/n=\infty$, then for
every real-valued function $f$ there exists a homeomorphism $h$
such that $\widehat{f\circ h}(k)=O(\varepsilon(|k|)/|k|+
1/|k|^{3/2})$. It is well-known that if $g$ is continuous and
$\widehat{g}(k)=o(1/|k|)$, then $g\in U(\mathbb T)$; thus
Sahakian's result implies the Bohr--P\'al theorem.

At the same time, in general, it is impossible to attain the
condition $\sum_{k\in\mathbb Z}|\widehat{f\circ h}(k)|<\infty$.
There exists a real-valued continuous function $f$ such that
whenever $h$ is a homeomorphism, the superposition $f\circ h$
does not belong to the space $A(\mathbb T)$ of absolutely
convergent Fourier series. This result, that provides a
solution to Lusin's rearrangement problem, is obtained by A. M.
Olevski\v{\i} [14] (see also [10]).

6. The Sahakian theorem, which guarantees  a very fast decrease
of the Fourier coefficients of the superposition $f\circ h$,
shows in particular, that for an appropriate homeomorphism we
have $f\circ h(k)=O(1/|k|)$ and $f\circ h\in W_2^{1/2}(\mathbb
T)$. At the same time, the Sahakian theorem does not imply even
the continuity of $\widetilde{f\circ h}$. It is easy to verify
that whenever $\{\varepsilon(n), n=1, 2, \ldots\}$ is a
sequence satisfying the assumptions of the Sahakian theorem,
one can find a real-valued continuous $g$ such that
$\widehat{g}(k)=O(\varepsilon(|k|)/|k|)$, but
$\widetilde{g}\notin L^\infty(\mathbb T)$. Indeed, consider the
series $\sum_{n=1}^\infty\varepsilon(n)n^{-1}\sin nt$. Since
$\varepsilon(n)$ monotonically decreases to zero, this series
converges uniformly. Let $g(t)$ be its sum. Then
$\widetilde{g}(t)\sim\sum_{n=1}^\infty-\varepsilon(n)n^{-1}\cos
nt$. Assuming $\widetilde{g}$ is in $L^\infty(\mathbb T)$, it
follows that the Fej\'er sums $\sigma_N(\widetilde{g})$ satisfy
$\|\sigma_N(\widetilde{g})\|_{L^\infty(\mathbb T)}=O(1)$. At
the same time
$$
\sigma_N(\widetilde{g})(t)=\sum_{n=1}^N\frac{-\varepsilon(n)}{n}\bigg(1-\frac{n}{N}\bigg)\cos nt,
$$
whence
$$
\|\sigma_N(\widetilde{g})\|_{L^\infty(\mathbb T)}=
\sum_{n=1}^N\frac{\varepsilon(n)}{n}\bigg(1-\frac{n}{N}\bigg).
\eqno(5)
$$
Since $\sum_{n=1}^\infty\varepsilon(n)/n=\infty$, the
right-hand side of (5) is unbounded.

\quad

I am grateful to A. B. Aleksandrov and V. Ya. Eiderman for a
useful discussion.

\begin{center}
\textbf{References}
\end{center}

\flushleft
\begin{enumerate}

\item H. Bohr, \"Uber einen Satz von J. P\'al, \emph{Acta
    Sci. Math. (Szeged)} \textbf{7} (1935), 129--135.

\item J. P\'al, Sur les transformations des fonctions qui
    font converger leurs s\'eries de Fourier, \emph{C. R.
    Acad. Sci Paris} \textbf{158} (1914), 101--103.

\item R. Salem, On a theorem of Bohr and P\'al, \emph{Bull.
    Amer. Math. Soc.}, \textbf{50} : 8 (1944), 579--580.

\item N. K. Bary, \emph{A treatise on trigonometric
        series}, Vols. I, II, Pergamon Press, Oxford 1964.

\item W. Jurkat, D. Waterman, Conjugate functions and the
    Bohr--P\'al theorem, \emph{Complex Variables} \textbf{12}
    (1989), 67--70.

\item G. Goffman, T. Nishiura, D. Waterman,
    \emph{Homeomorphisms in Analysis}, Mathematical Surveys
    and Monographs v. 54, Amer. Math. Soc., 1997.

\item  N. S. Stylianopoulos, E. Wegert, A uniform estimate
    for the modulus of continuity of starlike mappings,
    \emph{Ann. Univ. Mariae Curie-Sklodowska}, Sect. A.,
    \textbf{56} : 10 (2002), 97–103.

\item  A. Di Bucchianico, Banach algebras, logarithms, and
    polynomials of convolution type, \emph{Journal of
    Mathematical Analysis and Applications}, \textbf{156} : 1
    (1991), 253--273.

\item  J.-P. Kahane,  Y. Katznelson, S\'eries de Fourier des
    fonctions born\'ee, \emph{Studies in Pure Math.}, in
    Memory of Paul Tur\'an, Budapest, 1983, pp. 395--410.
    (Preprint, Orsay, 1978.)

\item  A. M. Olevskii, Modifications of functions and Fourier
    series, \emph{Uspekhi Mat. Nauk}, \textbf{40} : 3(243)
    (1985), 157–193; English transl.: \emph{Russian Math.
    Surveys}, \textbf{40} : 3 (1985), 181--224.

\item  V. V. Lebedev, Change of variable and the rapidity of
    decrease of Fourier coefficients, \emph{Mat. Sbornik},
    \textbf{181} : 8 (1990), 1099--1113; English transl.:
    \emph{Math. USSR-Sb.}, \textbf{70} : 2 (1991), 541--555.
    English transl. corrected by the author is available at:
    http://arxiv.org/abs/1508.06673

\item  V. Lebedev, The Bohr--P\'al theorem and the Sobolev
    space $W_2^{1/2}$, \emph{Studia Mathematica},
    \textbf{231} : 1 (2015), 73--81.

\item  A. A. Saakjan, Integral moduli of smoothness and the
    Fourier coefficients of the composition of functions,
    \emph{Mat. Sb. (N.S.)}, \textbf{110}(152):4(12) (1979),
    597–608; English transl: \emph{Mathematics of the
    USSR-Sbornik}, \textbf{38} : 4 (1981), 549--561.
	
\item  A. M. Olevskii, Change of variable and absolute
    convergence of Fourier series, \emph{Dokl. Akad. Nauk
    SSSR}, \textbf{256} : 2 (1981), 284--288; English
    transl.: \emph{Soviet Math. Dokl.}, \textbf{23} (1981),
76--79.

\end{enumerate}

\quad

\quad

National Research University Higher School of Economics

e-mail: \emph{lebedevhome@gmail.com}

\end{document}